\documentclass[12pt,reqno,smallextended]{amsart}
\usepackage{amsmath,amsfonts,amssymb}
\usepackage{t1enc}
 \numberwithin{equation}{section}

\begin{document}


\theoremstyle{plain}
\newtheorem*{theorem}{Theorem}
\newtheorem*{corollary}{Corollary}
\newtheorem{thm}{Theorem}[section]
\newtheorem{prop}[thm]{Proposition}
\newtheorem{lemma}[thm]{Lemma}
\newtheorem{sublem}[thm]{Sublemma}
\newtheorem{cor}[thm]{Corollary}
\newtheorem{dfn}[thm]{Definition}
\newtheorem{rem}[thm]{Remark}
\newtheorem*{que}{Question}
\newtheorem{example}[thm]{Example}


\newcommand{\lab}{\label}
\newcommand{\labeq}[1]{  \be \label{#1}  }
\newcommand{\labea}[1]{  \bea \label{#1}  }
\newcommand{\ben}{\begin{enumerate}}
\newcommand{\een}{\end{enumerate}}
\newcommand{\bm}{\boldmath}
\newcommand{\Bm}{\Boldmath}
\newcommand{\bea}{\begin{eqnarray}}
\newcommand{\ba}{\begin{array}}
\newcommand{\bean}{\begin{eqnarray*}}
\newcommand{\ea}{\end{array}}
\newcommand{\eea}{\end{eqnarray}}
\newcommand{\eean}{\end{eqnarray*}}
\newcommand{\beq}{\begin{equation}}
\newcommand{\eeq}{\end{equation}}
\newcommand{\bthm}{\begin{thm}}
\newcommand{\ethm}{\end{thm}}
\newcommand{\blem}{\begin{lem}}
\newcommand{\elem}{\end{lem}}
\newcommand{\bprop}{\begin{prop}}
\newcommand{\eprop}{\end{prop}}
\newcommand{\bcor}{\begin{cor}}
\newcommand{\ecor}{\end{cor}}
\newcommand{\brem}{\begin{rem}}
\newcommand{\erem}{\end{rem}}
\newcommand{\bdfn}{\begin{dfn}}
\newcommand{\edfn}{\end{dfn}}
\newcommand{\lb}{\linebreak}
\newcommand{\nlb}{\nolinebreak}
\newcommand{\nl}{\newline}
\newcommand{\hs}{\hspace}
\newcommand{\vs}{\vspace}
\alph{enumii} \roman{enumiii}
\def\endpf{\qed}
\def\ms{\medskip}
\def\endpf{\qed}
\def\ms{\medskip}
\def\gb{g_\beta}
\def\vep{\varepsilon}
\def\N{I\!\!N}                \def\Z{Z\!\!\!Z}      \def\R{I\!\!R}
\def\C{C\!\!\!\!I}            \def\T{T\!\!\!\!T}      \def\oc{\overline \C}
\def\Q{Q\!\!\!\!|}
\def\1{1\!\!1}
\def\pl{\wp_\Lambda}
\def\and{\text{ and }}        \def\for{\text{ for }}
\def\tif{\text{ if }}         \def\then{\text{ then }}
\def\Cap{\text{Cap}}          \def\Con{\text{Con}}
\def\Comp{\text{{\rm Comp}}}  \def\diam{\text{\rm {diam}}}
\def\dist{\text{{\rm dist}}}        \def\Dist{\text{{\rm Dist}}}   \def\Crit{{\rm Crit}}
\def\Sing{{\rm Sing}}          \def\Tr{{\rm Tr}}
\def\Fin{{\cal F}in}
\def\F{{\cal F}}
\def\h{{\text h}}
\def\hmu{\h_\mu}           \def\htop{{\text h_{\text{top}}}}
\def\H{\text{{\rm H}}}     \def\im{\text{{\rm Im}}} \def\re{\text{{\rm re}}}
\def\HD{\text{{\rm HD}}}   \def\DD{\text{{\rm DD}}}
\def\BD{\text{{\rm BD}}}         \def\PD{\text{PD}}
\def\HypDim{\text{{\rm HypDim}}}
\def\Int{\text{{\rm Int}}} \def\ep{\text{e}}
\def\CD{\text{CD}}         \def\P{\text{{\rm P}}}      \def\Id{\text{Id}}
\def\Ker{\text{{\rm Ker}}} \def\rank{\text{{\rm rank}}}\def\Mod{{\rm Mod}}
\def\HHH{\HD}              \def\const{{\rm const}}     \def\sh{\sharp}
\def\abs{\prec}
\def\A{\Cal A}             \def\Ba{\Cal B}             \def\Ca{\Cal C}
\def\Ha{\Cal H}
           \def\M{\Cal M}              \def\Pa{\Cal P}
\def\U{\Cal U}             \def\V{\Cal V}
\def\W{\Cal W}
\def\a{\alpha}                \def\b{\beta}             \def\d{\delta}
\def\De{\Delta}               \def\e{\epsilon}          \def\f{\phi}
\def\g{\gamma}                \def\Ga{\Gamma}
\def\La{\Lambda}              \def\om{\omega}           \def\Om{\Omega}
\def\Sg{\Sigma}               \def\sg{\sigma}
\def\Th{\Theta}               \def\th{\theta}           \def\vth{\vartheta}
\def\ka{\kappa}
\def\Ka{\Kappa}
\def\bi{\bigcap}              \def\bu{\bigcup}
\def\({\bigl(}                \def\){\bigr)}
\def\lt{\left}                \def\rt{\right}
\def\bv{\bivee}
\def\ld{\ldots}               \def\bd{\partial}         \def\^{\tilde}
\def\endpf{\hfill{$\endpfsuit$}}\def\Pr{\prod}
\def\tm{\widetilde{\mu}}
\def\tn{\widetilde{\nu}}
\def\es{\emptyset}            \def\sms{\setminus}
\def\sbt{\subset}             \def\spt{\supset}
\def\gek{\succeq}             \def\lek{\preceq}
\def\eqv{\Leftrightarrow}     \def\llr{\Longleftrightarrow}
\def\lr{\Longrightarrow}      \def\imp{\Rightarrow}
\def\comp{\asymp}
\def\upto{\nearrow}           \def\downto{\searrow}
\def\sp{\medskip}      \def\fr{\noindent}        \def\nl{\newline}
\def\ov{\overline}            \def\un{\underline}
\def\ds{\oplus}               \def\PC(F){{\rm PC}(F)}
\def\ess{{\rm ess}}
\def\ni{\noindent}
\def\om{\omega}
\def\bt{{\bf t}}
\def\Bu{{\bf u}}
\def\tr{t}
\def\bo{{\bf 0}}
\def\nut{\nu_\lla^{\scriptscriptstyle 1/2}}
\def\arg{\text{arg}}
\def\Arg{\text{Arg}}
\def\re{{\rm Re}}
\def\gr{\nabla}
\def\endpf{{${\mathbin{\hbox{\vrule height 6.8pt depth0pt width6.8pt  }}}$}}
\def\Fa{\cal F}
\def\Gal{\cal G}
\def\1{1\!\!1}
\def\D{{I\!\!D}}
   \def\PC(f){{\rm PC}(f)}
\def\La{\Lambda}
\def\vol{\text{\rm {vol}}}
\def\dc{\overline{\mathbb{C}}}
\def\oc{\mathbb{C}}
\def\int{{\rm int}}
\def\diam{\text{\rm {diam}}}   \def\vol{\text{\rm {vol}}}
\def\dist{\text{{\rm dist}}}        \def\Dist{\text{{\rm Dist}}}   \def\Crit{{\rm Crit}}
\def\dim{\text{{\rm dim}}}


\title[Hausdorff dimension of  elliptic functions]
{Hausdorff dimension of elliptic functions with critical values
approaching infinity}

\thanks{2000{\em Mathematics Subject Classification.} Primary 37F35.
Secondary 37F10, 30D05.\\The author is partially supported by PW -
grant No 504G/1120/0077/000 and MNiSW - grant No NN 201 607640}


\author{Piotr  Ga{\l}\k{a}zka}
\maketitle

\smallskip
{\footnotesize
 \centerline{Faculty of Mathematics and Information Sciences}
   \centerline{Warsaw University of Technology}
   \centerline{Warsaw 00-661, Poland}
   \centerline{{\it tel.: }+48 22 234 75 40, {\it fax.: }+48 22 625 74 60}
   \centerline{{\it e-mail: }P.Galazka@mini.pw.edu.pl}
}

\begin{abstract}
We   consider the {\it escaping parameters} in the family $\b \pl$,
i.e. these parameters $\b$ for
 which  the orbits of  critical values  of $\b \pl$   approach
 infinity, where  $\pl$ is the Weierstrass  function. Unlike to the exponential map
 the considered functions are ergodic. They admit  a non-atomic, $\sigma$-finite, ergodic,
conservative and invariant measure $\mu$ absolutely continuous with
respect  to the  Lebesgue measure. Under additional assumptions on
the $\pl$-function we estimate from below the Hausdorff dimension of
the set of escaping parameters in the family $\b \pl$, and compare
it with the Hausdorff dimension of escaping set in dynamical space,
proving a similarity between parameter plane and dynamical space.
\end{abstract}

\maketitle

\section{Introduction}

 In the series of papers J.~Hawkins and L.~Koss \cite{HK3,
HK1,HK2} described dynamics of Weierstrass functions. Ergodic theory
of non-recurrent elliptic functions was developed by J.~Kotus and
M.~Urba\'nski in \cite{KU1,KU2,KU3}. Recently, in \cite{HKK} there were given examples of all possible behaviours of non-recurrent elliptic functions (called in that paper
critically tame functions). These include the map with critical
values approaching infinity.  The aim of this paper is to show that
the escaping parameters form a considerably big set.\\

  Let $f\colon\oc\to\dc$ be a transcendental meromorphic function. For
$n\in\mathbb{N}$, denote by $f^n$ the $n$-th iterate of $f$. The
\textit{Fatou set} $F(f)$ of $f$ is the set of points $z\in\oc$ such
that all iterates $f^n(z)$ are well-defined and the family
$\{f^n\}_{n\in\mathbb{N}}$ is a normal family in some neighbourhood
of $z$. The complement of $F(f)$ in $\dc$ is called the
\textit{Julia set} of $f$.  P.~Dom\'\i nguez in
 \cite{D} proved  that for transcendental meromorphic functions with poles {\it the escaping set} $$I(f)=\{ z \in {\mathbb C}: \,
\lim_{n\to\infty} f^n(z)=\infty\}$$ is
 not  empty and $ J(f)=\partial{I(f)}$.
 Later P.~Rippon and G.~Stallard \cite{RS} showed that if, additionally a function  $f$ is in
 the  Eremenko-Lyubich  class $\mathcal B$, then $I(f) \sbt J(f)$.
 It means that $\Int I(f)=\es $. Recently, several authors \cite{BKS,BK,BRS,RS1,RS2} have studied
properties of the escaping set for entire and meromorphic functions.
 In \cite{K1} the  Hausdorff dimension of $I(f)$ was estimated  from  below for some  class of meromorphic functions.
 In particular, this  can be applied to elliptic functions of the form
 $\gb=\b\pl$, $ \b \in \mathbb C \sms\{0\} $, where $\pl$ is the Weierstrass elliptic  function. As a corollary we obtain  that the  Hausdorff dimension
 $\dim_H (I(\gb))\geq 4/3$.  On the other
hand, Bergweiler, Kotus and Urba\'nski proved in \cite{BK,KU1} that an upper bound on
$\dim_H (I(\gb))$ is the same as a lower bound,
 so
\begin{equation*}
 \dim_H (I(\gb))=\frac{4}{3}.
 \end{equation*}
 In this paper, we additionally  assume  that a lattice of
$\pl$-function  is triangular and the  critical values of $\pl$
 are the poles.  As a counterpart of escaping set $I(\gb)$ we consider the
set of escaping parameters in the family $\gb$, i.e.
$${\mathcal E}=\{ \b  \in
{\mathbb C}\sms\{0\} :\,\, \lim_{n\to\infty} \gb^n(c_i)= \infty,\,\,
i=1,2,3 \},$$ where $c_i$ is a  critical point of $\pl$.
 For these
maps the Julia set is the whole plane $\overline{\mathbb C}$. In
this paper we construct a
 collection of Cantor subsets of ${\mathcal E}$ with prescribed rate of growth
 and estimate from below their Hausdorff dimension.
 The main result is the following theorem.

\begin{theorem}\label{glowne2}
For any  one-parameter family of functions
  $\gb(z)=\b\pl(z),$
where  $\b \in \mathbb C \sms \{0\}$, $\La=[ \lambda_1, e^{2\pi i/3}
\lambda_1 ]$ is a triangular lattice such that all critical values
of $\pl$ are the poles, the Hausdorff dimension of the set of
escaping parameters $\mathcal E$ is greater or equal to 4/3.
\end{theorem}

The paper is organized as follows. In section 2 we give the
background definitions and results for studying elliptic
functions, in particular the $\pl$-Weierstrass function. We also
summarize  metric properties  of maps in $\mathcal E$. In sections 3
and 4 we show how one can find escaping parameters. In the last
section we  estimate from below $\dim_H(\mathcal E)$.

\

\section{General preliminaries}

We start with recalling the definition and basic properties of
elliptic functions. For $\lambda_1,\lambda_2\in\oc\setminus\{0\}$
such that $\im(\lambda_1/\lambda_2)\neq 0$ a lattice $\La\sbt\oc$ is
defined by $$\La=[\lambda_1,\lambda_2]=\{l\lambda_1+m\lambda_2,\,
l,m\in\mathbb{Z}\}.$$
\begin{dfn}
An elliptic function is a meromorphic function $f\colon\oc\to\dc$
which is periodic with respect to a lattice $\La$, i.e.
$
f(z)=f(z+l\lambda_1+m\lambda_2)
$
for all $z\in\oc$ and $l,m\in\mathbb{Z}$.
\end{dfn}

  We denote by $b_{l,m}=l \lambda_1+m \lambda_2$,\,
$l,m\in\mathbb{Z}$, lattice points of $\Lambda$ and by
$$\mathcal{R}=\{t_1 \lambda_1+t_2 \lambda_2; \,\, 0\leq t_1,t_2 < 1\}$$ a
fundamental parallelogram of $\La$. For a non-constant elliptic
function and a given $w\in \overline{\mathbb{C}}$ the number of
solutions of the equation $f(z)=w$ in $\mathcal{R}$ equals the sum
of multiplicities of the poles in a fundamental parallelogram. Since
the derivative of an elliptic function is also an elliptic function
periodic with respect to the same lattice, then each elliptic
function has infinitely many critical points but only finitely many
critical values. Due to periodicity elliptic functions do not have
asymptotic values. Thus they belong to the class $\mathcal{S}$.

A special case of an elliptic function is the Weierstrass elliptic
function defined by
$$
\pl (z)=\frac{1}{z^2}+\sum_{w\in\La\setminus\{0\}} \left(
\frac{1}{(z-w)^2}-\frac{1}{w^2}\right)
$$
for all $z\in\oc$ and every lattice $\La$. It is well-known that
$\pl$ is periodic with respect to $\La$ and has order 2. The
derivative of the Weierstrass function is also an elliptic function
periodic with respect  to $\La$  and is defined by
$$
\pl' (z)=-2 \sum_{w\in\La} \frac{1}{(z-w)^3}.
$$
The Weierstrass elliptic function and its derivative are related by
the differential equation
\begin{equation}\label{rr}
\left(\pl' (z)\right)^2=4 \left(\pl (z)\right)^3-g_2 \pl (z)- g_3,
\end{equation}
where $g_2= g_2(\La)= 60 \sum_{w\in\La\setminus\{0\}}
\frac{1}{w^4},\, g_3= g_3(\La)= 140 \sum_{w\in\La\setminus\{0\}}
\frac{1}{w^6}$. The numbers $g_2(\La),g_3(\La)$ are invariants of
the lattice $\La$ in the following sense, if $g_i(\La)=g_i(\La'),\,
i=2,3$, then $\La=\La'$. Moreover, for any $g_2,g_3$ such that
$g_2^3-27g_3^2\neq 0$ there is a lattice $\La$ with invariants
$g_2,g_3$. For any lattice $\La$ the Weierstrass function $\pl$
satisfies the property of homogeneity, i.e.
\begin{equation}\label{Wjedn}
\wp_{\a\La}(\a z)=\frac{1}{\a^2}\pl(z)
\end{equation}
for every $\a\in\oc\setminus\{0\}$. The Weierstrass function has
poles of order 2 at lattice points and its derivative has poles of
order 3. In the fundamental parallelogram  the map $\pl$ has three critical points
which we denote by
$$
c_1=\frac{\lambda_1}{2},\ c_2=\frac{\lambda_2}{2},\
c_3=\frac{\lambda_1+\lambda_2}{2}.
$$
We use the symbols $e_i=\pl(c_i),\, i=1,2,3$ to denote the critical
values of $\pl$. They are related to each other with the equations
\begin{equation}\label{wkrytw}
e_1+e_2+e_3=0,\ e_1e_3+e_2e_3+e_1e_2=-\frac{g_2}{4},\
e_1e_2e_3=\frac{g_3}{4}.
\end{equation}
We consider only the Weierstrass functions periodic with respect to
triangular lattices, i.e. lattices $\La=[\lambda_1,\lambda_2]$ such
that $\lambda_2=e^{2\pi i/3}\lambda_1$. In other words a lattice is
triangular if $\La=e^{2\pi i/3}\La$. For triangular lattices $g_2=0$
and the critical values of $\pl$ are the cube roots of $g_3/4$.
Moreover, \eqref{rr} and \eqref{wkrytw} imply that the critical
value $e_3$ is a non-zero real number and $e_1,e_2$ are given by the
formulas $e_1=e^{4\pi i/3}e_3,\, e_2=e^{2\pi i/3}e_3$. The
iterations of the critical values turn out to have the same
property, i.e. $\pl^n(e_1)=e^{4\pi i/3}\pl^n(e_3),\,
\pl^n(e_2)=e^{2\pi i/3}\pl^n(e_3),\, n\geq 1$. It is a consequence
of an invariance of a triangular lattice with respect to the
rotation $z\mapsto e^{2\pi i/3}z$ and the homogeneity of $\pl$ given
in \eqref{Wjedn} (see \cite{HK1} for details).

\

 We additionally assume that all the
critical values of the Weierstrass function $\pl$ are poles.  The
example of a family of such lattices was given by J. Hawkins i L.
Koss in \cite{HK1}.
\begin{example}\label{ex0}
Let $\Omega=[\omega_1,\omega_2]$ be a lattice with invariants
$g_2=0,g_3=4$. It is a triangular lattice for which $e_1=e^{4\pi
i/3},e_2=e^{2\pi i/3},e_3=1$. Let $\gamma_1=\sqrt[3]{\frac{e^{4\pi
i/3}\omega_1^2}{m}}$, where $m$ is an odd negative number and
$\gamma_2=\frac{\gamma_1 \om_2}{\omega_1}$. Then, the lattice
$\Gamma=[\gamma_1,\gamma_2]$ is triangular and all the critical
values of $\pl$ are poles.
\end{example}

\

Now, we describe ergodic properties of so-called critically tame
elliptic functions studied by J. Kotus and M. Urba\'nski in
\cite{KU3}. We start with some definitions and notations.

\begin{dfn}
Let $f\colon\oc\to\dc$ be an elliptic function and $z\in\oc$ such
that all iterates $f^n(z),n\in\mathbb{N}$ are well-defined. A point
$w\in\dc$ is called an $\om-$\textit{limit} point of $z$ for $f$, if
there is a sequence of natural numbers $n_k\to\infty$ such that
$$
\lim_{k\to\infty} \dist_s (f^{n_k}(z),w)=0,
$$
where $\dist_s$ denotes spherical metric in $\dc$. The
$\om-$\textit{limit} set of $z$ is a set of all $\om-$\textit{limit}
points of $z$ and we denote it by $\om(z)$.
\end{dfn}

\begin{dfn}\label{pc}
Let $D$ be a domain in $\oc$ and $g\colon D\to\oc$ an analytic map.
Set $z\in\oc,\ r>0$. We denote by $U(z,g^{-1},r)$ the connected
component of $g^{-1}(B(g(z),r))$containing $z$. Suppose that $c\in
Crit(g)$. Then, there exist $r=r(g,c)>0$ and $K=K(g,c)\geq 1$ such
that
$$
\frac{1}{K}|z-c|^{p}\leq |g(z)-g(c)| \leq K|z-c|^{p}
$$
and
$$
\frac{1}{K}|z-c|^{p-1}\leq |g'(z)| \leq K|z-c|^{p-1}
$$
for all $z\in U(c,g^{-1},r)$ and some natural $p=p(g,c)$, and also
such that
$$
g(U(c,g^{-1},r))=B(g(c),r).
$$
The number $p$ is called the order of $g$ at the critical point $c$
and is denoted by $p_c$. The number $p_c-1$ is the multiplicity of
the zero of $g'$ at $c$.
\end{dfn}

Denote  by   ${\mathcal P}_n(f)$, $ n\geq 1$,  the set   of prepoles
  of  order $n$ of $f$, i.e. $$ {\mathcal P}_n(f)=\{z\in\oc\colon
f^n(z)=\infty\}.$$ In particular, ${\mathcal P}_1(f)$ is the set of
poles of $f$.

\begin{dfn}
Suppose that $f\colon\oc\to\dc$ is an elliptic function and $b\in
{\mathcal P}_1(f)$. Let $\eta_b$ denote the multiplicity of the pole
$b$. We define
$$
q:=\sup \{\eta_b\colon b\in {\mathcal P}_1(f)\}=\max \{\eta_b\colon
b\in {\mathcal P}_1(f)\cap \mathcal{R}\}.
$$
\end{dfn}

\noindent Denote by $\Crit (f)$ the set of critical points of $f$,
i.e.
$$ \Crit (f)=\{z\in\oc\colon f'(z)=0\}. $$
Let $\Crit_b (f)$ be the set of all prepole critical points, i.e.
$$
\Crit_b (f)=\Crit (f)\cap \bigcup_{n\in\mathbb{N}} {\mathcal
P}_n(f).
$$
Moreover, we define the set of all critical points of $f$ which
trajectories approach infinity, i.e.
$$
\Crit_\infty (f)=\{c\in\Crit(f)\colon \lim_{n\to\infty}
f^n(c)=\infty\}.
$$
Note that ${\mathcal P}_n(f)=f^{-1}({\mathcal P}_{n-1}(f))$ for all
$n\geq 2$ and ${\mathcal P}_n(f)\sbt J(f)$. For every
$c\in\Crit_b(f)$ there is a unique $n\in\mathbb{N}$ such that $c\in
{\mathcal P}_n(f)$. For all $c\in\Crit_\infty(f)$ and every $R>0$
there exists natural $N$ such that for all $n\geq N:
|f^{n+1}(c)|>R$. This inequality is equivalent to the fact that
$f^n(c)$ lies close to a unique pole $b_n$. That implies that for
all $c\in\Crit_\infty(f)$ one can define a sequence of poles $b_n$
close to the iterates of $f$.

\begin{dfn}
Let $f\colon\oc\to\dc$ be an elliptic function. For $c\in
\Crit_\infty(f)$ we define
$$
q_c:=\limsup_{n\to\infty} \eta_{b_n},
$$
where the sequence $\{b_n\}_{n\geq 1}$ was defined above. Moreover,
let
$$
l_\infty=\max \{p_c q_c\colon c\in\Crit_\infty (f)\},
$$
where $p_c$ is as in Definition \ref{pc}.
\end{dfn}

\begin{dfn}\label{crittame}
Let $f\colon \oc\to\dc$ be an elliptic function and $c\in\Crit(f)$.
We say that $f$ is critically tame if the following conditions are
satisfied:
\begin{itemize}
    \item[(a)] if $c\in F(f)$, then there exists an attracting or parabolic cycle of period $p$, $S=\{z_0,f(z_0),\ldots,f^{p-1}(z_0)\}$ such that $\om(c)=S$.
    \item [(b)] if $c\in J(f)$, then one of the following holds:
\begin{itemize}
    \item [(i)] $\om(c)$ is a compact subset of $\oc$ such that $c\notin\om(c)$,
    \item [(ii)] $c\in\Crit_b (f)$,
    \item [(iii)] $c\in\Crit_\infty (f)$ and
    $$ \dim_H (J(f))>\frac{2 l_\infty}{l_\infty +1}. $$
\end{itemize}
\end{itemize}
\end{dfn}

\

Denote by $\Tr(f)\index{Tr(f)@$Tr(f)$}\sbt J(f)$ the set of all
transitive points of\index{set of  transitive points} $f$, that is
the set of points in $J(f)$ such that their forward trajectories are
dense in $J(f)$.

\

We quote two results from \cite{KU3}, which became an inspiration
for studying the escaping parameters $\mathcal E$. Below  a
conformal measure $m$ is defined by means of the spherical metric.

\begin{prop}\label{tmaincm}  Suppose that $f$ is a critically tame
elliptic function, denote $h=\dim_H (J(f))$. Then there exist:
\begin{itemize}
\item [a)]  a unique atomless $h$-conformal measure $m$ for $f\colon J(f)\sms\{\infty\} \to J(f)$, $m$ is ergodic, conservative and $m(\Tr(f))=1$.
\item [b)]  a non-atomic, $\sigma$-finite, ergodic, conservative and
invariant measure $\mu$ for $f$, equivalent to the  measure $m$.
Additionally, $\mu$ is unique up to a multiplicative constant and is
supported on $J(f)$. \end{itemize}
\end{prop}

The next proposition gives  sufficient conditions for an elliptic
function $f$ to satisfy the conditions given in Definition
\ref{crittame}.

\begin{prop}\label{Hwardost}
If every critical point $c$ of $f$ is such that $c\in \Crit_b (f)$
or $c\in \Crit_\infty (f)$, then $J(f)=\dc$ and $f$ is critically
tame.
\end{prop}

Proposition~\ref{tmaincm} and Proposition~\ref{Hwardost} imply that
the elliptic functions considered  in the next sections are ergodic
with respect to the Riemann measure $m$. This is in contrast with
Lyubich's result \cite{Lu2} which says  that $e^z$ is not ergodic
with respect to the Lebesgue measure.  The escaping parameters in
the exponential family $f_\lambda(z)=\lambda e^z$, \,$ \lambda \in
\mathbb C\setminus\{0\}$, were also studied by Urba\'nski and Zdunik
in \cite{UZ}. Under the assumption that the forward trajectory of
$0$ grows exponentially fast (this includes the case $ \lambda >
1/e$), they showed that $\om(z)=\{f^n_\lambda(0):\, n \geq 0\}\cup
\{\infty\}$ for a.e. $z \in J(f_\lambda)=\ov{\mathbb C}$. Later
Hemke \cite{H} proved that these maps are non-recurrent. His results
cover the fast escaping parameters in the tangent family
$f_\lambda(z)=\lambda \tan(z), \lambda \in \mathbb C\setminus\{0\}$,
for which again he proved that $\om(z)=\{f^n_\lambda(\pm \lambda
i):\, n \geq 0\}\cup \{\infty\}$ for a.e. $z \in
J(f_\lambda)=\ov{\mathbb C}$. In all the cases the existence of  a
non-atomic, $\sigma$-finite, ergodic, conservative and invariant
measure $\mu$ for $f$, absolutely continuous with respect to  the
Lebesgue measure  follows from \cite{KU} or
Proposition~\ref{tmaincm}.

\

At the end of this section we recall a definition of distortion. Let
$U$ be an open subset of $\oc$, $f\colon U\to\oc$ be a conformal
map, then its distortion  is defined as
$$ L(f,U):=\frac{\sup_{z\in U}|f'(z)|}{\inf_{z\in U}|f'(z)|}.$$
For conformal maps we have
\begin{equation}\label{dystodw}
 L(f,U)=L(f^{-1},f(U)).
\end{equation}
To prove a lower bound on $\dim_H(\mathcal E)$ we use the following
theorem proved by C.~McMullen in \cite{McM}.
\begin{prop}\label{mcmullen} For each $n\in\mathbb{N}$, let 
$\mathcal{A}_n$ be a finite collection of disjoint compact
subsets of $\mathbb{R}^d$, each of which has positive $d$-dimensional
Lebesgue measure. Define
$$
\mathcal{U}_n=\bigcup_{A_n\in\mathcal{A}_n} A_n,\ \
A=\bigcap_{n=1}^\infty \mathcal{U}_n.
$$
Suppose that for each $A_n\in\mathcal{A}_n$ there is $A_{n+1}
\in\mathcal{A}_{n+1}$ and a unique $A_{n-1} \in\mathcal{A}_{n-1}$
such that $A_{n+1}\sbt A_n\sbt A_{n-1}$. If $\Delta_n,d_n$ are such that, for each $A_n\in\mathcal{A}_n$,
\begin{align*}
& \frac{\vol (\mathcal{U}_{n+1}\cap A_n)}{\vol (A_n)}\geq \Delta_n>0, \\
& \diam (A_n)\leq d_n <1, \\
& d_n \stackrel{n\to\infty}{\longrightarrow} 0,
\end{align*}
then
$$
\dim_H (A)\geq d-\lim_{n\to\infty} \sum_{j=1}^n \frac{|\log
\Delta_j|}{|\log d_n|}.
$$
\end{prop}

\

\section{The escaping parameters}

Unlike to the exponential or tangent family we do not know any
examples of  Weierstrass functions with critical values approaching
infinity. In this section, we recall from \cite{HKK} how one can find
the elliptic functions with critical values eventually mapped onto
poles (Lemma~\ref{konstr2}) and the maps with critical values
escaping to infinity (Lemma~\ref{konstr3}).

We  consider one-parameter family of functions
  $$\gb(z)=\b\wp_\Lambda(z),$$
where  $\b \in \mathbb C \sms \{0\}$, $\La=[ \lambda_1, e^{2\pi i/3}
\lambda_1 ]$ is a triangular lattice such that all critical values
of $\wp_\La$ are the poles. These lattices were constructed  in
\cite{HK1} (see also Example~\ref{ex0}). The functions under
consideration $\gb$ are periodic and their critical points are the
same as for the Weierstrass function $\pl$.  It was shown in
\cite{HKK} that the critical orbits of $\gb$ behave symmetrically,
i.e.
\begin{equation}\label{symetria}
g^n_\b(c_2)=\gamma^2 g^n_\b(c_1),\,\,\, g^n_\b(c_3)=\gamma
g^n_\b(c_1)
\end{equation}
for all $n \in \mathbb N$, where $\gamma=e^{2\pi i/3}$.  So we can
take only one of them. Let it be the trajectory of the critical
value $\gb(c_1)$. Denote $B_\rho(\infty):=\{z\in\dc\colon
|z|>\rho\},\rho>0$. To prove the next lemma we consider the
auxiliary functions $h_n(\b)=\gb^n(c_1), n \in \mathbb N$. It will
appear in the proof of the next lemma that these functions are
defined outside a countable set of parameters.

\begin{lemma}\label{konstr2}
Let $\La$ be a triangular lattice such that all critical values of
$\pl$ are the poles. For every $r>0$ and each $n\geq 2$, there is
$\b\in B(1,r)$,  such that $\gb^n (c_1)=\infty$.
\end{lemma}
\fr{\sl Proof.} Consider the function $h_1$ defined before, i.e.
$h_1\colon B(1,r)\to\oc$, $h_1(\b)=\gb(c_1)$, where $0<r<1/2$. By the
assumption, $h_1(1)=g_1(c_1)=\pl(c_1)$ is a pole of $\pl$. Now we
define $h_2\colon B(1,r)\to\dc$ by the formula $h_2(\b)=\gb^2(c_1)$.
Denote by $\mathcal{P}(h_2)$ the set of its poles. Since
$h_2(1)=g_1^2(c_1)=\pl^2(c_1)=\infty$, then $1\in \mathcal{P}(h_2)$.
Thus, the theorem is true for $n=2$. We can take $r$ so small that
$1$ is a unique pole of $h_2$ in $B(1,r)$. Actually, let $\b\in
B(1,r)\setminus\{1\}$ be a pole of $h_2$. Thus,
$h_2(\b)=\gb^2(c_1)=\b\pl(\b\pl(c_1))=\infty$, so
$\pl(\b\pl(c_1))=\infty$, which implies $\b\pl(c_1)\in\La$. However
$\pl(c_1)\in\La$, so taking $r$ small enough we have
$\b\pl(c_1)\notin\La$ for $\b \in B(1,r)\sms \{1\}$. Then, $h_2$ is
a non-constant meromorphic function. Since $1$ is a pole (of order
2) of the function $h_2$, then we can take $R_2\geq 2^2$ such that
$B_{R_2}(\infty)\sbt h_2(B(1,r))$. The set $B_{R_2}(\infty)$
contains infinitely many lattice points $b_{l,m}^{(2)}$ of $\La$ and
each of them (being a pole of $\pl$) is the image of some parameter
$\b_{l,m}^{(2)}\in B(1,r)\setminus\{1\}$ under $h_2$. Choose one of
$\b_{l,m}^{(2)}$ and denote it, for simplicity, by $\b_2$. We denote
the corresponding pole by $b_2$. We have constructed the map
$g_{\b_2}$, such that the orbit of the critical point $c_1$ is the
following
$$
c_1 \mapsto g_{\b_2}(c_1) \mapsto g^2_{\b_2}(c_1)=b_2 \mapsto
g^3_{\b_2}(c_1)=\infty,
$$
where $g_{\b_2}(c_1)$ is close to (but not equal to) the critical
value $\pl(c_1)$ and $g^2_{\b_2}(c_1)\in B_{R_2}(\infty)$. Let
$r_1:=r$. Take $0<r_2<r_1/2$ so small that $\overline{B(\b_2,r_2)}\sbt B(1,r)\setminus \mathcal{P}(h_2)$ and $h_3(B(\b_2,r_2))\sbt B_{R_2}(\infty)$, where $h_3(\b)=\gb^3(c_1)$. Restricting $h_3$ to $B(\b_2,r_2)$, we take $R_3\geq 2R_2\geq 2^3$ such that $B_{R_3}(\infty)\sbt h_3(B(\b_2,r_2))$. Each lattice point $b_{l,m}^{(3)}\in B_{R_3}(\infty)$ is the image of some parameter $\b_{l,m}^{(3)}\in B(\b_2,r_2)\setminus\{\b_2\}$. Note that this proves the existence of a parameter $\b_3$ such that
$$
c_1 \mapsto g_{\b_3}(c_1)\approx \pl(c_1) \mapsto g^2_{\b_3}(c_1)\approx b_2 \mapsto g^2_{\b_3}(c_1)= b_3 \mapsto
g^4_{\b_3}(c_1)=\infty,
$$
where none of the $\approx$ are equality and $b_i\in \La\cap B_{R_i}(\infty)$ with $R_i\geq 2^i,\ i=2,3$.
Now, by induction we define a map with the property that the
critical point is a prepole of order $n\geq 4$.   Fix $n\geq 4$
and suppose for all $k<n$ we have constructed the maps
$$
h_k\colon B(1,r)\setminus \bigcup_{1<i<k} \mathcal{P}(h_i) \to \dc
$$
by the formulas $h_k(\b)=\gb^k(c_1)$, where $\mathcal{P}(h_i)$ is
the set of poles of $h_i$. We define a map
$$
h_n\colon B(1,r)\setminus \bigcup_{1<k<n} \mathcal{P}(h_k) \to \dc
$$
such that $h_n(\b)=\gb^n(c_1)$. The set $\bigcup_{1<k<n}
\mathcal{P}(h_k)$ is a set of essential singularities of $h_n$. In
its complement the map $h_n$ is meromorphic, denote by
$\mathcal{P}(h_n)$ its set of poles. Set a pole
$\b_{n-1}\in\mathcal{P}(h_n)$. The equality
$h_n(\b_{n-1})=g^n_{\b_{n-1}}(c_1)=\infty$ implies that there is a small enough
constant $0<r_{n-1}<r_{n-2}/2$ such that $\overline{B(\b_{n-1},r_{n-1})}\sbt B(\b_{n-2},r_{n-2})\setminus \bigcup_{1<k<n} \mathcal{P}(h_k)$ and $h_n(B(\b_{n-1},r_{n-1}))\sbt B_{R_{n-1}}(\infty)$. Now, we can take $R_n\geq 2R_{n-1}\geq 2^n$ such that $B_{R_n}(\infty)\sbt h_n(B(\b_{n-1},r_{n-1}))$. Next, we choose
one of the lattice points of $\La$ from $B_{R_n}(\infty)$ and denote
it by $b_n$. We know that $b_n$ is the image of some parameter
$\b_n\in B(\b_{n-1},r_{n-1})\setminus\{\b_{n-1}\}$, i.e. $b_n=h_n(\b_n)=g^n_{\b_n}
(c_1)$. The orbit of the critical point $c_1$ for the map $g_{\b_n}$
is the following
$$
c_1 \mapsto g_{\b_n}(c_1)\approx \pl(c_1) \mapsto
g^2_{\b_n}(c_1)\approx b_2 \mapsto \ldots \mapsto
g^n_{\b_n}(c_1)=b_n \mapsto g^{n+1}_{\b_n}(c_1)=\infty,
$$
where $g^i_{\b_n}(c_1)\in B_{R_i}(\infty),i=1,\ldots ,n$. This
finishes the proof.
\qed

\

\begin{lemma}\label{konstr3}
Let $\La$ be a triangular lattice such that all critical values of
$\pl$ are the poles. Then, for every $r>0 $  there is a parameter
$\b\in B(1,r)$ such that $\lim_{n\to\infty}\gb^n (c_i)=\infty,\
i=1,2,3$.
\end{lemma}
\fr{\sl Proof.} We show that $\lim_{n\to\infty}\gb^n (c_1)=\infty$.
The 'symmetry' of the critical orbits given in (\ref{symetria})
implies the lemma is true for $c_2$ and $c_3$. By
Lemma~\ref{konstr2}, there is a sequence of parameters
$\{\b_n\}_{n\geq 2}$ such that
$$
\left|g^n_{\b_n}(c_1)\right|>R_n\ \mbox{ and }\
g^{n+1}_{\b_n}(c_1)=\infty,
$$
where $R_n\geq 2^n$ and a decreasing sequence of balls
$B(\b_n,r_n)\sbt B(1,r_1)\setminus\bigcup_{1<k<n}\mathcal{P}(h_k)$
such that $r_n< 2^{-n}$. Since $r_n\to 0$, then there is the
parameter $\b=\bigcap_{n\geq 2} \overline{B(\b_n,r_n)}$. By the
construction from the proof of Lemma~\ref{konstr2}, $\b$ is an
accumulation point of the set $\bigcup_{n>1}\mathcal{P}(h_n)$. The
iterates of the critical point under $\gb$ satisfy the conditions
$\left|g^n_{\b}(c_1)\right|>R_n\geq 2^n$ for all $n\geq 2$. Hence,
$\lim_{n\to\infty} R_n=\infty$, which implies
$\lim_{n\to\infty}\gb^n (c_1)=\infty$. \qed

\

\section{Escaping parameters with prescribed rate of growth of
critical orbits}

In this section, we construct a collection  of  subsets of $
\mathcal E$ with prescribed rate of growth of the critical orbits of
$\gb$. We  fix a function $\pl$ such that
\begin{equation*}
 \La=[ \lambda_1, e^{2\pi i/3} \lambda_1 ]
 \end{equation*}
is a triangular lattice and all critical values of $\pl$  are the
poles. These lattices were constructed  in \cite{HK1} (see also
Example~\ref{ex0}). We consider one-parameter family of functions
\begin{equation}\label{family}
\gb(z)=\b\wp_\Lambda(z), \,\, \b \in B(1,r) \,\, \mbox{ for} \,\,
0<r <1/2.
\end{equation}
The functions $\gb$ are periodic and their critical points are the
same as for the Weierstrass function $\pl$. It follows from
(\ref{symetria}) that the critical orbits of $\gb$ behave
symmetrically, i.e.
$$
g^n_\b(c_2)=\gamma^2 g^n_\b(c_1),\,\,\, g^n_\b(c_3)=\gamma
g^n_\b(c_1)
$$
for all $n \in \mathbb N$, where $\gamma=e^{2\pi i/3}$. Since $\pl$ is periodic, there exist
  a constant
\begin{equation*}
 0 < \vep_0 < \min\{1,  |\lambda_1|/3\}
 \end{equation*}
 and holomorphic functions $G,H$
such that for each pole $b_{l,m}\in\La$
$$\begin{aligned}
\pl(z) = &\frac{a_{-2}}{(z-b_{l,m})^2}+ \frac{a_{-1}}{z-b_{l,m}}+
\sum_{k=0}^\infty a_k(z-b_{l,m})^k
      = : \frac{G(z)}{(z-b_{l,m})^2}\\
 \pl'(z) = & \frac{b_{-3}}{(z-b_{l,m})^3}+
\frac{b_{-2}}{(z-b_{l,m})^2}+ \frac{b_{-1}}{z-b_{l,m}} +
\sum_{k=0}^\infty b_k(z-b_{l,m})^k
      =:  \frac{H(z)}{(z-b_{l,m})^3}\end{aligned}$$
for all $z\in B(b_{l,m},\vep_0)$, where $G(b_{l,m})=a_{-2}\neq 0$,
$H(b_{l,m})=b_{-3}\neq 0$. Shrinking $\varepsilon_0$, if necessary,
we may assume that $G(z)\neq 0$ and  $H(z)\neq 0$ for $z \in
B(b_{l,m},\varepsilon_0)$. The periodicity of $\pl$ implies
that there exist universal constants $K_1,K_2>0$ such that
\begin{equation*}
K_1^{-1}\leq |G(z)|\leq K_1, \quad  K_2^{-1}\leq |H(z)|\leq K_2
\end{equation*}
on all balls $B(b_{l,m},\vep_0)$. Hence,
\begin{equation*}
\frac{K_1^{-1}}{|z-b_{l,m}|^2}\leq |\pl(z)|=\left|\frac{G(z)}{(z-b_{l,m})^2}\right| \leq \frac{K_1}{|z-b_{l,m}|^2}
\end{equation*}
and
\begin{equation*}
\frac{K_2^{-1}}{|z-b_{l,m}|^3}\leq |\pl'(z)|=\left|\frac{H(z)}{(z-b_{l,m})^3}\right| \leq
\frac{K_2}{|z-b_{l,m}|^3}
\end{equation*}
for all $l,m\in\mathbb{Z}$ and $z\in B(b_{l,m},\vep_0)$. For every
$\b\in B(1,r)$, where $r$ is defined in \eqref{family} and for all
$z\in B(b_{l,m},\vep_0),l,m\in\mathbb{Z}$, we have
\beq\label{szacgb}
  \frac{C_1^{-1}}{|z-b_{l,m}|^2}\leq
\left|g_\b(z)\right|=|\b \pl(z)|\leq \frac{C_1}{|z-b_{l,m}|^2} \eeq
and \beq\label{szacgbder} \frac{C_2^{-1}}{|z-b_{l,m}|^3}\leq
\left|g_\b'(z)\right|=|\b \pl'(z)|\leq \frac{C_2}{|z-b_{l,m}|^3}
\eeq where $C_1=2K_1,C_2=2K_2$. Moreover, shrinking $\vep_0,r$ if
necessary, we can choose constants $M_1,M_2,\, 0<M_2-M_1<\pi/4$ such
that \beq\label{staleM} M_1\leq \arg(\b G(z))\leq M_2 \eeq for all
$\b\in B(1,r)$ and $z\in B(b_{l,m},\vep_0),l,m\in\mathbb{Z}$. We
recall from Section 3 that
$$ h_1\colon B(1,r)\to  \mathbb C, \,\, h_1(\b)= \gb(c_1), $$
where $c_1  $ is a  critical point of $\pl$. We choose
$\vep >0$  such that the following conditions are  simultaneously satisfied
\begin{equation}\label{epsilon2}
\begin{aligned}
& \vep <\min\{\vep_0, |\pl(c_1)|/3\},\\
& B(\pl(c_1), \vep)\sbt h_1(B(1, r)),\\
& \pl \mbox{ is one-to-one on each of the segments defined in \eqref{wycinek}.}
\end{aligned}
\end{equation}
Let
\begin{equation}\label{wycinek}
 U(z_0,\vep):=\{z\in\oc\colon -\frac{3\pi}{8}\leq \Arg (z-z_0)\leq\frac{3\pi}{8},\
|z-z_0|\leq\vep\}, \end{equation}
where $z_0\in\La$ and $\vep$ is defined above. Next, we take  $R_1 >0$  such that
$$ U(\pl(c_1), \vep) \subset P(0, R_1, 2R_1):=\{ z\in \mathbb
C : R_1< |z| < 2R_1 \}.$$
Using \eqref{szacgb} and \eqref{staleM}, we get
\begin{equation*}
\begin{aligned}
& \{z\in\dc\colon |z|\geq\frac{C_1}{\vep^2},\ -\frac{3\pi}{4}+M_2\leq \arg z\leq\frac{3\pi}{4}+M_1 \}  \sbt \gb(U(b_{l,m},\vep))\\ & \sbt \{z\in\dc\colon |z|\geq\frac{C_1^{-1}}{\vep^2},\ -\frac{3\pi}{4}+M_1\leq \arg z\leq\frac{3\pi}{4}+M_2 \}
\end{aligned}
\end{equation*}
for all $l,m\in\mathbb{Z}$. Since $0<M_2-M_1<\pi/4$, there exists $\phi\in\mathbb{R}$ such that
\beq\label{zawieraniekonstr}
\{z\in\dc\colon |z|\geq\frac{C_1}{\vep^2},\ \phi-\frac{\pi}{8}\leq \arg z\leq\phi+\frac{9\pi}{8} \}  \sbt \gb(U(b_{l,m},\vep)).
\eeq
We choose $R_2$ such that
\beq\label{stalaR2}
R_2>\frac{C_1}{(1-\alpha)\vep^2},
\eeq
where $\alpha=\sin(\pi/8)=\sqrt{2-\sqrt{2}}/2$.
Thus, it  follows from \eqref{zawieraniekonstr} that
\begin{equation}\label{pokrycie}
 \{z\in \mathbb C:  |z| > R_2,\ \phi\leq \arg z\leq\phi+\pi\}\sbt \gb(U(b_{l,m}, \vep))
\end{equation}
 for all the poles $b_{l,m}$. Let $a_1=R_2/R_1> \frac{C_1}{(1-\alpha)\vep^2 R_1}.$
Now, we define a constant
\begin{equation}\label{stala-a0}
 a_0= \max\left\{ 2, a_1, \frac{1}{R_1},\frac{3C_1^{3/2}}{C_2 R_1},\frac{6^4 C_1^6}{C_2^4
R_1^5},  \left(\frac{4 \vep (1+r)C_1^{3/2}}{C_2
R_1^{5/2}}\right)^{2/3},\frac{\sqrt{C_1}}{\sqrt[3]{C_2} \sqrt{R_1}}
\right\}.
\end{equation}
Fix
\begin{equation*}
a > a_0
\end{equation*} and consider a sequence or radii
\begin{equation*}
R_n:=a^{n-1}R_1, \quad  n \geq 2.
\end{equation*}
Let
\begin{equation*}
 P(0, R_n, 2R_n):=\{z\in\oc\colon R_n<|z|<2R_n \}, \quad n \geq 2
\end{equation*} and
\begin{equation}\label{polpierscien}
P^+(0,R_n,2R_n):=\{z\in\oc\colon R_n<|z|<2R_n,\ \phi<\arg
z<\phi+\pi\}, \quad n \geq 2.
\end{equation}
The condition $a>a_0\geq 2$ guarantees that the annuli $P(0, R_n,
2R_n)$ are pairwise disjoint. Recall that in the previous section we
defined the auxiliary functions $h_n(\b)=g_\b^n(c_1)$, $n \in
\mathbb N$.

\

\begin{dfn}\label{konstr4}
We define the following family of sets
\begin{align*}
\mathcal{A}_0(a) & = \{A_0=B(1,r)\}, \\
\mathcal{A}_1(a) & = \{A_1=h_1^{-1}(U(\pl(c_1),\vep))\sbt A_0\}, \\
\mathcal{A}_2(a) & = \{A_2\sbt A_1\ |\ \exists
b_{l,m}^{(2)}\in\La\colon  U(b_{l,m}^{(2)},\vep)\sbt
P^+(0,R_2,2R_2),\ A_2=
h_2^{-1}(U(b_{l,m}^{(2)},\vep))\}, \\
\ldots & \\
\mathcal{A}_n(a) & = \{A_n\sbt A_{n-1}\ |\ \exists
b_{l,m}^{(n)}\in\Lambda\colon  U(b_{l,m}^{(n)},\vep)\sbt
P^+(0,R_n,2R_n),\
A_n=h_n^{-1}(U(b_{l,m}^{(n)},\vep))\}, \\
\ldots &
\end{align*}
where $h_n^{-1}(U(b_{l,m}^{(n)},\vep))$ denotes a component of the
preimage under the map $h_n^{-1}$.  Let
$$ \mathcal{U}_n(a)=\bigcup_{A_n\in\mathcal{A}_n(a)} A_n,\ \ 
A(a)=\bigcap_{n=1}^\infty \mathcal{U}_n(a).$$
\end{dfn}

\

\bprop\label{nonempty} For each $n \in \mathbb N$ the set ${\mathcal
A}_n(a) $ defined above is non-empty. \eprop

\fr{\sl  Proof.} In the previous section, we
showed that the function $h_2$ has the pole at
$\b=1=h_1^{-1}(\pl(c_1))\in\partial A_1$. Thus,
$\mathcal{A}_1(a)\neq\emptyset$. Since $h_1(A_1)=U(\pl(c_1),\vep)$, it follows from \eqref{pokrycie} that
$$h_2(A_1)=\{\gb(h_1(\b))|\b\in A_1\} \supset P^+(0, R_2, 2R_2).$$
Take a pole $b^{(2)}_{l,m} \in \La\cap P^+(0, R_2, 2R_2)$ with $U(b_{l,m}^{(2)},\vep)\sbt P^+(0,R_2,2R_2)$. Since
$h_2(A_1)\supset P^+(0, R_2, 2R_2)$  there exists $\b^{(2)}_{l,m}\in A_1$
such that $h_2(\b^{(2)}_{l,m})=b_{l,m}^{(2)}$. Thus, the set $\mathcal
A_2(a)$  is non-empty. Now, we fix $n\geq 3$ and suppose that
$\mathcal{A}_{n-1}(a)\neq\emptyset$. We will show that $\mathcal A_n(a)
\neq \es$.  Since $h_{n-1}(A_{n-1})=U(b_{l,m}^{(n-1)},\vep)$ for some $b_{l,m}^{(n-1)}\in\La\cap P^+(0,R_{n-1},2R_{n-1})$, it follows from \eqref{pokrycie} that
$$h_{n}(A_{n-1})=\{\gb(h_{n-1}(\b))|\b\in A_{n-1}\} \supset  P^+(0, R_{n}, 2R_{n}),$$
as $R_n=a^{n-2}R_2$ and $a>a_0\geq 2$ in view of \eqref{stala-a0}.
Choosing $\b_{l,m}^{(n)}\in A_{n-1}$ such that $h_{n}(\b_{l,m}^{(n)})=b_{l,m}^{(n)} \in
\La\cap P^+(0, R_{n}, 2R_{n})$ and $U(b_{l,m}^{(n)},\vep)\sbt P^+(0,R_n,2R_n)$, we obtain that $\mathcal A_n(a)\neq\emptyset $. By induction, the lemma is true for all $n\in\mathbb{N}$. \qed

\

\bthm\label{Cantor} Let $\gb$ be the family of  maps defined above
and let  $a_0$ be a constant given  in \eqref{stala-a0}. Then, for
every $a > a_0$
 there is a Cantor subset $A(a)$ of
$\mathcal E$ such that
$$\dim_H(A(a))\geq \frac{4}{3} - \frac{6 \log 2}{\log a}.$$
\ethm \bcor For $ a  \to + \infty$  we  have $\dim_H(A(a)))\nearrow
\frac{4}{3}$ and $\dim_H(\mathcal E) \geq \frac{4}{3}.$ \ecor

\section{The proofs}

In this section  we prove Theorem~\ref{Cantor}. We  fix $a > a_0$
and consider the sets $\mathcal A_n(a)$, $ n \geq 1$, defined in
Definition~\ref{konstr4}. We  drop the parameter $a$ and  keep
notation from the last section.\\

The first  two lemmas  are devoted to estimates of the derivatives
$h'_n$, $ n \geq 2$.

\begin{lemma}\label{hd1}
Let $A_n \in \mathcal A_n$, $n \geq 2$. Then, for every $\b\in A_n$
$$
h_n'(\b)=\frac{1}{\b}\prod_{k=1}^{n-1} \gb'(\gb^k(c_1))\left[
\gb(c_1)+\sum_{k=2}^n \frac{\gb^k (c_1)}{\prod_{i=1}^{k-1}
\gb'(\gb^i(c_1))} \right].
$$
\end{lemma}
\fr{\sl Proof.} Let  $n=2$. Then,
$$
\begin{aligned}
    h_1(\b)&=\gb(c_1)=\b \pl(c_1), \\
    h_2(\b)&=\gb^2(c_1)=\b \pl(\b \pl(c_1)), \\
    h_2'(\b)&= \pl(\b \pl(c_1)) + \b \pl'(\b \pl(c_1))\pl(c_1) = \frac{\gb^2(c_1)}{\b}+\frac{\gb'(\b\pl(c_1))\gb(c_1)}{\b}\\
    &= \frac{1}{\b} \gb'(\b\pl(c_1)) \left[ \gb(c_1)+ \frac{\gb^2(c_1)}{\gb'(\gb(c_1))}\right]= \frac{1}{\b} \gb'(\gb(c_1)) \left[ \gb(c_1)+ \frac{\gb^2(c_1)}{\gb'(\gb(c_1))}\right].
    \end{aligned}
    $$
Suppose that the lemma is true for some $n\geq 2$. We show that it
is true for $n+1$.
  $$  \begin{aligned}
    h_{n+1}(\b) & =\b\pl(h_n(\b)), \\
    h_{n+1}'(\b)&= \pl(h_n(\b)) + \b\pl'(h_n(\b))\cdot h_n'(\b)\\
    =& \frac{\gb^{n+1}(c_1)}{\b} + \gb'(\gb^n(c_1))\cdot \frac{1}{\b}\cdot \prod_{k=1}^{n-1} \gb'(\gb^k(c_1))\cdot \left[ \gb(c_1)+
    \sum_{k=2}^n \frac{\gb^k (c_1)}{\prod_{i=1}^{k-1} \gb'(\gb^i(c_1))} \right]\\
    =& \frac{\gb^{n+1}(c_1)}{\b} + \frac{1}{\b}\cdot \prod_{k=1}^{n} \gb'(\gb^k(c_1))\cdot \left[ \gb(c_1)+
    \sum_{k=2}^n \frac{\gb^k (c_1)}{\prod_{i=1}^{k-1} \gb'(\gb^i(c_1))} \right]\\
    =& \frac{1}{\b}\cdot \prod_{k=1}^{n} \gb'(\gb^k(c_1))\cdot \left[ \gb(c_1)+
    \sum_{k=2}^n \frac{\gb^k (c_1)}{\prod_{i=1}^{k-1} \gb'(\gb^i(c_1))} + \frac{\gb^{n+1} (c_1)}{\prod_{i=1}^{n} \gb'(\gb^i(c_1))} \right]\\
    =& \frac{1}{\b}\cdot \prod_{k=1}^{n} \gb'(\gb^k(c_1))\cdot \left[ \gb(c_1)+\sum_{k=2}^{n+1} \frac{\gb^k (c_1)}{\prod_{i=1}^{k-1} \gb'(\gb^i(c_1))} \right].
    \end{aligned}
$$
\qed

\

 \fr We recall from the previous section (see \eqref{szacgb},\eqref{szacgbder}) that there are
universal constants $C_1,\, C_2>0$ such that
$$
  \frac{C_1^{-1}}{|z-b_{l,m}|^2}\leq
\left|g_\b(z)\right|\leq \frac{C_1}{|z-b_{l,m}|^2},\ \
\frac{C_2^{-1}}{|z-b_{l,m}|^3}\leq \left|g_\b'(z)\right|\leq
\frac{C_2}{|z-b_{l,m}|^3}
$$ for all $l,m\in\mathbb{Z}$, every $z\in
B(b_{l,m},\vep)$ and all $\b\in B(1,r)$. To simplify the formulas in
the following part of the paper we write
\begin{equation}\label{uzyteczne}
|\gb(z)| \asymp \frac{C_1}{|z-b_{l,m}|^2},\ \ |\gb'(z)| \asymp
\frac{C_2}{|z-b_{l,m}|^3}.
\end{equation}
\ni Note that if $\b\in\mathcal{U}_n,\, n\geq 2$ and $z=\gb^j(c_1)$ with $j\in\{1,2,\ldots n-1\}$ we have $\gb(z)=\gb^{j+1}(c_1)=h_{j+1}(\b)\in U(b_{l,m}^{(j+1)},\vep)\sbt P^+(0,R_{j+1},2R_{j+1})$ and moreover, using \eqref{uzyteczne},
\begin{equation}\label{asymp1}
R_{j+1} \leq |\gb(z)| \asymp \frac{C_1}{|z-b_{l,m}|^2} \leq 2R_{j+1}
\end{equation}
for some $b_{l,m}\in \La\cap P^+(0,R_{j},2R_{j})$.
The inequality \eqref{asymp1} implies that
$$
\frac{C_1}{2R_{j+1}} \leq |z-b_{l,m}|^2 \leq \frac{C_1}{R_{j+1}},
$$
which is equivalent to
$$
\left(\frac{C_1}{2R_{j+1}}\right)^{3/2} \leq \left|z-b_{l,m}\right|^3
\leq \left(\frac{C_1}{R_{j+1}}\right)^{3/2}.
$$
Then,
$$
\frac{C_2}{\left(\frac{C_1}{R_{j+1}}\right)^{3/2}} \leq |\gb'(z)| \asymp
\frac{C_2}{\left|z-b_{l,m}\right|^3} \leq
\frac{C_2}{\left(\frac{C_1}{2R_{j+1}}\right)^{3/2}}
$$
or, equivalently,
\begin{equation}\label{nierpochH}
\frac{C_2 R_{j+1}^{3/2}}{C_1^{3/2}} \leq |\gb'(z)| \leq \frac{2^{3/2}
C_2 R_{j+1}^{3/2}}{C_1^{3/2}}
\end{equation}
for $\b\in\mathcal{U}_n,\, n\geq 2$ and $z=\gb^j(c_1)$ with $j\in\{1,2,\ldots n-1\}$.

\

\begin{lemma}\label{hd2}
Let $A_n \in \mathcal A_n$, $n \geq 2$. Then,  for every $\b\in A_n$
$$
\frac{1}{2(1+r)} \left(\frac{C_2}{C_1^{3/2}}\right)^{n-1}
a^\frac{3n(n-1)}{4} R_1^\frac{3n-1}{2} \leq \left|h_n'(\b)\right|
\leq \frac{5}{2(1-r)} \left(\frac{2^{3/2}
C_2}{C_1^{3/2}}\right)^{n-1} a^\frac{3n(n-1)}{4}
R_1^\frac{3n-1}{2}.
$$
\end{lemma}

\fr{\sl Proof.} In Lemma \ref{hd1}, we proved that
$$
h_n'(\b)=\frac{1}{\b}\prod_{k=1}^{n-1} \gb'(\gb^k(c_1))\left[
\gb(c_1)+\sum_{k=2}^n \frac{\gb^k (c_1)}{\prod_{i=1}^{k-1}
\gb'(\gb^i(c_1))} \right]
$$
for all $n\geq 2$ and every $\b\in A_n$. First, we estimate the
product $\prod_{k=1}^{n-1}\gb'(\gb^k(c_1))$. Observe that $$
\gb(\gb^k(c_1))=\gb^{k+1}(c_1)=h_{k+1}(\b),\ k=1,2,\ldots ,n-1.$$
The functions $h_2,\ldots h_n$ are well-defined for $\b\in A_n$,
because $A_n\sbt A_k,\, k=2,\ldots, n$. Since $h_{k+1}(\b)\in
P(0,R_{k+1},2R_{k+1})$, then using \eqref{nierpochH}, we get
\begin{align*}
\left|\prod_{k=1}^{n-1}\gb'(\gb^k(c_1)\right| & \leq \frac{2^{3/2} C_2 R_2^{3/2}}{C_1^{3/2}}\cdot\ldots\cdot \frac{2^{3/2} C_2 R_n^{3/2}}{C_1^{3/2}}\\& = \left(\frac{2^{3/2} C_2}{C_1^{3/2}}\right)^{n-1} (aR_1)^{3/2}\cdot\ldots\cdot (a^{n-1}R_1)^{3/2} \\
& = \left(\frac{2^{3/2} C_2}{C_1^{3/2}}\right)^{n-1}
a^\frac{3n(n-1)}{4} R_1^\frac{3(n-1)}{2}.
\end{align*}
Analogously, we get the estimate from below
$$
\left|\prod_{k=1}^{n-1}\gb'(\gb^k(c_1)\right| \geq
\left(\frac{C_2}{C_1^{3/2}}\right)^{n-1} a^\frac{3n(n-1)}{4}
R_1^\frac{3(n-1)}{2}.
$$
Finally,
\begin{equation}\label{hd2p1}
\left(\frac{C_2}{C_1^{3/2}}\right)^{n-1} a^\frac{3n(n-1)}{4}
R_1^\frac{3(n-1)}{2} \leq
\left|\prod_{k=1}^{n-1}\gb'(\gb^k(c_1)\right| \leq
\left(\frac{2^{3/2} C_2}{C_1^{3/2}}\right)^{n-1} a^\frac{3n(n-1)}{4}
R_1^\frac{3(n-1)}{2}.
\end{equation}
Now, using \eqref{hd2p1}, we estimate the sum $\sum_{k=2}^n
\frac{\gb^k (c_1)}{\prod_{i=1}^{k-1} \gb'(\gb^i(c_1))}$.
$$
\begin{aligned}
 & \left|\sum_{k=2}^n \frac{\gb^k(c_1)}{\prod_{i=1}^{k-1} \gb'(\gb^i(c_1))}\right| \leq  \sum_{k=2}^n
\left| \frac{\gb^k(c_1)}{\prod_{i=1}^{k-1} \gb'(\gb^i(c_1))} \right|
\leq
\sum_{k=2}^n \frac{2R_k}{\left(\frac{C_2}{C_1^{3/2}}\right)^{k-1} a^\frac{3k(k-1)}{4} R_1^\frac{3(k-1)}{2}} \\
 & = \sum_{k=2}^n
\frac{2a^{k-1}R_1}{\left(\frac{C_2}{C_1^{3/2}}\right)^{k-1}
a^\frac{3k(k-1)}{4} R_1^\frac{3(k-1)}{2}}=
\sum_{k=2}^n \frac{2}{\left(\frac{C_2}{C_1^{3/2}}\right)^{k-1} a^\frac{(k-1)(3k-4)}{4} R_1^\frac{3k-5}{2}} \\
 & = \frac{2 C_1^{3/2}}{C_2 \sqrt[4]{a R_1}} \sum_{k=2}^n
\left(\frac{C_1^{3/2}}{C_2}\right)^{k-2}
\frac{1}{a^\frac{3k^2-7k+3}{4} R_1^\frac{6k-11}{4}}.
\end{aligned}
$$
Since $a>a_0\geq 2$ and $3k^2-7k+3 \geq 6k-11 $ for
$k=2,3,\ldots$, then
$$
\sum_{k=2}^n \left(\frac{C_1^{3/2}}{C_2}\right)^{k-2}
\frac{1}{a^\frac{3k^2-7k+3}{4} R_1^\frac{6k-11}{4}} \leq
\sum_{k=2}^n \left(\frac{C_1^{3/2}}{C_2}\right)^{k-2} \frac{1}{(a
R_1)^\frac{6k-11}{4}}.
$$
Using the inequality $(6k-11)/4\geq k-2$ for $k\geq 3/2$ and the fact that $a>a_0\geq \max\{\frac{1}{R_1},\frac{3C_1^{3/2}}{C_2 R_1}\}$, we get
$$
\sum_{k=2}^n \left(\frac{C_1^{3/2}}{C_2}\right)^{k-2} \frac{1}{(a
R_1)^\frac{6k-11}{4}} \leq \sum_{k=2}^n \left(\frac{C_1^{3/2}}{C_2 a
R_1}\right)^{k-2} \leq \sum_{k=2}^\infty \left(\frac{C_1^{3/2}}{C_2
a R_1}\right)^{k-2} = \frac{1}{1-\frac{C_1^{3/2}}{C_2 a R_1}} <
\frac{3}{2}.
$$

\fr Hence,
\begin{equation}\label{hd2p2}
\left|\sum_{k=2}^n \frac{\gb^k (c_1)}{\prod_{i=1}^{k-1}
\gb'(\gb^i(c_1))}\right| \leq \frac{2 C_1^{3/2}}{C_2 \sqrt[4]{a
R_1}}\cdot \frac{3}{2} = \frac{3 C_1^{3/2}}{C_2 \sqrt[4]{a R_1}}
\leq \frac{R_1}{2},
\end{equation}
because $a>a_0\geq \frac{6^4 C_1^6}{C_2^4 R_1^5}$. Using
\eqref{hd2p2}, we get
\begin{equation}\label{hd2p3}
\frac{R_1}{2} = R_1 - \frac{R_1}{2} \leq \left|
\gb(c_1)+\sum_{k=2}^n \frac{\gb^k (c_1)}{\prod_{i=1}^{k-1}
\gb'(\gb^i(c_1))} \right| \leq 2R_1 + \frac{R_1}{2}= \frac{5
R_1}{2}.
\end{equation}
Plugging \eqref{hd2p1}, \eqref{hd2p3} into the formula for $h'_n$
from Lemma \ref{hd1}, we obtain
$$
|h_n'(\b)| \leq \frac{1}{1-r} \left(\frac{2^{3/2}
C_2}{C_1^{3/2}}\right)^{n-1} a^\frac{3n(n-1)}{4}
R_1^\frac{3(n-1)}{2} \cdot \frac{5R_1}{2} = \frac{5}{2(1-r)}
\left(\frac{2^{3/2} C_2}{C_1^{3/2}}\right)^{n-1} a^\frac{3n(n-1)}{4}
R_1^\frac{3n-1}{2}
$$
and
$$
|h_n'(\b)| \geq \frac{1}{1+r}
\left(\frac{C_2}{C_1^{3/2}}\right)^{n-1} a^\frac{3n(n-1)}{4}
R_1^\frac{3(n-1)}{2} \cdot \frac{R_1}{2} = \frac{1}{2(1+r)}
\left(\frac{C_2}{C_1^{3/2}}\right)^{n-1} a^\frac{3n(n-1)}{4}
R_1^\frac{3n-1}{2}
$$
for $a>a_0$. Both estimates prove the lemma.
\qed

\

In Proposition \ref{nonempty}, we showed that each set
$\mathcal{A}_n$, defined in Definition \ref{konstr4}, is non-empty
and its elements, the sets $A_n$, contain on their boundary
parameters $\b_n$ such that $h_n(\b_n)\in \La\cap P(0,R_n,2R_n)$. In
the next part of this section, we estimate the  diameters of $A_n$
and the ratios $ \vol (\mathcal{U}_{n+1}\cap A_n)/\vol
(A_n)$. To do that  we should know that the functions $h_n$ are
conformal on $A_n \in \mathcal{A}_n $. Note that the maps $h_n,
n\geq 2$, are holomorphic outside a countable set of points and have
poles at $\b_{n-1}\in\partial A_{n-1}$.

\

\begin{lemma}\label{konstr6}
For each  $A_n\in{\mathcal A}_n$, $ n \geq 1$,  the map $h_n$ is
conformal on $A_n$.
\end{lemma}

\fr{\sl Proof.} The map $h_1$ is one-to-one and holomorphic
on $A_1$. By induction, we show that the maps $h_n,n\geq 2$ are conformal. Suppose that $h_n,n\geq 1$ is conformal on $A_n$, we prove that $h_{n+1}$ is conformal on $A_{n+1}\sbt A_n$. If $n=1$ then we take a segment
$$ U(b^{(1)}_{l,m},\vep)\sbt P(0,R_1,2R_1) $$
with $b^{(1)}_{l,m}=\pl(c_1)$ and if $n\geq 2$ we consider a segment
$$ U(b^{(n)}_{l,m},\vep)\sbt P^+(0,R_n,2R_n).$$
We know that $A_n=h_n^{-1}(U(b^{(n)}_{l,m},\vep)),n\geq 1$. Let
$b^{(n)}_{l,m}=b_n,\, \b_n=h_n^{-1}(b_n)\in\partial A_n$ and
$b^{(n+1)}_{l,m}=b_{n+1}$. If $U(b_{n+1},\vep)\sbt
P^+(0,R_{n+1},2R_{n+1})$, then
$h_{n+1}^{-1}(U(b_{n+1},\vep))=A_{n+1}\sbt A_n$.  We define a map
$\tilde{h}_{n+1}(\b)=\b_n \pl(h_n(\b))$. It follows from
\eqref{zawieraniekonstr} and \eqref{stalaR2} that
$$
\tilde{h}_{n+1}(A_n)\supset \{z\in\dc\colon |z|\geq\frac{C_1}{\vep^2},\ \phi-\frac{\pi}{8}\leq \arg z\leq\phi+\frac{9\pi}{8} \}
$$
and (shrinking $r$ if necessary) there is a set $\tilde{A}_{n+1}$ such that $A_{n+1}\sbt\tilde{A}_{n+1}\sbt A_n$ and
\begin{equation}\label{zawrozn}
\tilde{h}_{n+1}(\tilde{A}_{n+1})= \{z\in\oc\colon (1-\alpha)R_{n+1}<|z|<(2+\alpha)R_{n+1},\ \phi-\frac{\pi}{8}\leq \arg z\leq\phi+\frac{9\pi}{8} \}
\end{equation}
for some $\phi\in\mathbb{R}$ and $\alpha=\sqrt{2-\sqrt{2}}/2$. We show that $\tilde{h}_{n+1}$ is one-to-one on $\tilde{A}_{n+1}$. Take $\b',\b''\in \tilde{A}_{n+1}$ such that $\tilde{h}_{n+1}(\b')=\tilde{h}_{n+1}(\b'')$.
By definition of the map $\tilde{h}_{n+1}$, we have
$\pl(h_n(\b'))=\pl(h_n(\b''))$, where $h_n(\b'),h_n(\b'')\in h_n(\tilde{A}_{n+1})\sbt h_n(A_n)=
U(b_n,\vep)$. Since $\pl$ is one-to-one on $h_n(\tilde{A}_{n+1})$,
then $h_n(\b')=h_n(\b'')$ and this implies that $\b'=\b''$. This
follows from the injectivity of the map $h_n$. It follows from \eqref{polpierscien} and \eqref{zawrozn} that
$$
\begin{aligned}
& U(b_{n+1},\vep) =h_{n+1}(A_{n+1})\sbt P^+(0,R_{n+1},2R_{n+1})\\ & \sbt \{z\in\oc\colon (1-\alpha)R_{n+1}<|z|<(2+\alpha)R_{n+1},\ \phi-\frac{\pi}{8}\leq \arg z\leq\phi+\frac{9\pi}{8}\} = \tilde{h}_{n+1}(\tilde{A}_{n+1}).
\end{aligned}
$$
Moreover, for $\b\in\partial\tilde{A}_{n+1},\, 0<r<1/4-1/(2\alpha+4)$ we have
$$ \dist(\tilde{h}_{n+1}(\b),h_{n+1}(A_{n+1}))\geq \alpha R_{n+1}>2r|\pl(\zeta)|,$$
where $\zeta=h_n(\b)\in\partial h_n(\tilde{A}_{n+1})$.\\
We define auxiliary maps $H_{n+1}(\b)=h_{n+1}(\b)-w,\
\tilde{H}_{n+1}(\b)=\tilde{h}_{n+1}(\b)-w$ with $w\in h_{n+1}(A_{n+1})$. Thus, for $\b\in\partial \tilde{A}_{n+1}$ we have
$$
|{\tilde H}_{n+1}(\b)|=|\tilde{h}_{n+1}(\b)-w|\geq
\dist(\tilde{h}_{n+1}(\b),h_{n+1}(A_{n+1}))>2r
|\pl(\zeta)|$$ and
\begin{align*}
|H_{n+1}(\b)-{\tilde H}_{n+1}(\b)|= & |h_{n+1}(\b)-{\tilde
h}_{n+1}(\b)|=|\b\pl(\zeta)-\b_n\pl(\zeta)|=\\ &
|\b-\b_n||\pl(\zeta)|<2r |\pl(\zeta)|.
\end{align*}
Hence, $|{\tilde H}_{n+1}(\b)|>|H_{n+1}(\b)-{\tilde H}_{n+1}(\b)|$
on the set $\partial \tilde{A}_{n+1}$. Since the map $h_{n+1}$
is holomorphic on $\int A_n$, then the maps
$H_{n+1},\tilde{H}_{n+1}$ are holomorphic on $\int \tilde{A}_{n+1}$ and
continuous on  $\partial \tilde{A}_{n+1}$. Thus, the assumptions of
Rouch\'e Theorem are satisfied. It implies that $\tilde{H}_{n+1}$
and $H_{n+1}=\tilde{H}_{n+1}+H_{n+1}-\tilde{H}_{n+1}$ have the same
number of zeros on $\tilde{A}_{n+1}$, or, equivalently, the equations
$\tilde{h}_{n+1}(\b)=w$ and $h_{n+1}(\b)=w$ have the same number of
roots in $\tilde{A}_{n+1}$. Since the map $\tilde{h}_{n+1}$ is
one-to-one on $\tilde{A}_{n+1}$, then the former equation has a unique
root for a given $w$. Thus, the latter as well. This proves that
$h_{n+1}$ is one-to-one on $A_{n+1}$. The map
$h_{n+1}$ is holomorphic on $\int A_n$, then is
conformal on $A_{n+1}$.
\qed

\

\begin{lemma}\label{hd3}
Let $A_n\in \mathcal{A}_n$, $ n \geq 2$. Then
$$
L(h_n,A_n)\leq \frac{5(1+r)}{1-r}\cdot 2^\frac{3(n-1)}{2}.
$$
\end{lemma}\index{dystorsja}

\fr{\sl Proof.} Using the definition of distortion and Lemma
\ref{hd2}, for $a>a_0$ we get
$$
L(h_n,A_n)=\frac{\sup_{\b\in A_n} |h'_n(\b)|}{\inf_{\b\in A_n}
|h'_n(\b)|}\leq \frac{\frac{5}{2(1-r)} \left(\frac{2^{3/2}
C_2}{C_1^{3/2}}\right)^{n-1} a^\frac{3n(n-1)}{4}
R_1^\frac{3n-1}{2}}{\frac{1}{2(1+r)}
\left(\frac{C_2}{C_1^{3/2}}\right)^{n-1} a^\frac{3n(n-1)}{4}
R_1^\frac{3n-1}{2}} = \frac{5(1+r)}{1-r}\cdot 2^\frac{3(n-1)}{2}.
$$
\qed

\

\begin{lemma}\label{hd4}
For each  $A_n\in \mathcal{A}_n$, $n \geq 2$,
$$
\diam (A_n) \leq \frac{4 \vep
(1+r)}{\left(\frac{C_2}{C_1^{3/2}}\right)^{n-1} a^\frac{3n(n-1)}{4}
R_1^\frac{3n-1}{2}},
$$
where $\vep$ is as in \eqref{epsilon2}.
\end{lemma}

\fr{\sl Proof.} From Definition \ref{konstr4} we know that each set
of the form $h_n(A_n)$ is a segment of radius $\vep$, so $\diam
(h_n(A_n))\leq 2\vep$. Using Lemma \ref{hd2}, for $a>a_0$ we get
$$
\diam (A_n)\leq \frac{\diam (h_n(A_n))}{\inf_{\b\in A_n} |h'_n(\b)|}
\leq \frac{2\vep}{\frac{1}{2(1+r)}
\left(\frac{C_2}{C_1^{3/2}}\right)^{n-1} a^\frac{3n(n-1)}{4}
R_1^\frac{3n-1}{2}} = \frac{4 \vep
(1+r)}{\left(\frac{C_2}{C_1^{3/2}}\right)^{n-1} a^\frac{3n(n-1)}{4}
R_1^\frac{3n-1}{2}}.
$$
\qed

\

\brem Observe that $\diam (A_n)\to  0$ as $n \to \infty$, since
$a>a_0\geq 2$. This proves that the set $A$ from Definition
\ref{konstr4} is a Cantor set of parameters. \erem

\

\ni By Lemma \ref{hd4}, the numbers $d_n$
defined in Proposition~\ref{mcmullen} are equal to
\begin{equation}\label{srednica}
d_n = \frac{4 \vep (1+r)}{\left(\frac{C_2}{C_1^{3/2}}\right)^{n-1}
a^\frac{3n(n-1)}{4} R_1^\frac{3n-1}{2}},\ n\geq 2
\end{equation}
and $d_1= \diam (A_1) \leq 2r<1$ by (\ref{family}). We have
$$ d_2= \frac{4 \vep (1+r)}{\frac{C_2}{C_1^{3/2}}
a^{3/2} R_1^{5/2}} = \frac{4 \vep (1+r) C_1^{3/2}}{C_2 a^{3/2}
R_1^{5/2}}.$$
A straightforward calculation shows that the condition
$d_2<1$ is equivalent to
$$ a > \left(\frac{4 \vep (1+r)C_1^{3/2}}{C_2 R_1^{5/2}}\right)^{2/3}. $$
Using \eqref{srednica}, we get
$$ \frac{d_{n+1}}{d_n} = \frac{\frac{4
\vep (1+r)}{\left(\frac{C_2}{C_1^{3/2}}\right)^{n}
a^\frac{3n(n+1)}{4} R_1^\frac{3n+2}{2}}}{\frac{4 \vep
(1+r)}{\left(\frac{C_2}{C_1^{3/2}}\right)^{n-1} a^\frac{3n(n-1)}{4}
R_1^\frac{3n-1}{2}}} = \frac{C_1^{3/2}}{C_2 a^{3n/2} R_1^{3/2}}
$$
and
$$ \frac{d_3}{d_2} =  \frac{C_1^{3/2}}{C_2 a^{3} R_1^{3/2}} <1 \ \Longleftrightarrow\ a^3 > \frac{C_1^{3/2}}{C_2 R_1^{3/2}} \ \Longleftrightarrow\  a > \frac{\sqrt{C_1}}{\sqrt[3]{C_2} \sqrt{R_1}}.
$$
Since $a>a_0\geq \max\left\{1,\left(\frac{4 \vep (1+r)C_1^{3/2}}{C_2
R_1^{5/2}}\right)^{2/3},\frac{\sqrt{C_1}}{\sqrt[3]{C_2}
\sqrt{R_1}}\right\}$ and $d_{n+1}/d_n<d_3/d_2$ for $n\geq 3$, we get
$d_n<1, n=2,3,\ldots$ as required in Proposition~\ref{mcmullen}.

\

Next, we estimate from below the density of the sets
$\mathcal{U}_{n+1}\cap A_n$ in the set $A_n \in \mathcal{A}_n$ for
all $n\geq 1$.

\begin{lemma}\label{hd5}
There exists $M >0$ such that
$$ \frac{\vol (\mathcal{U}_{n+1}\cap
A_n)}{\vol (A_n)} \geq \frac{M}{2^{9n} R_{n+1}},
$$
for each  $A_n\in \mathcal{A}_n$, $n\geq 2$.
 Moreover,
$$
\frac{\vol (\mathcal{U}_2\cap A_1)}{\vol (A_1)} \geq \frac{M'}{R_2},
$$
for some $M'>0$.
\end{lemma}

\fr{\sl Proof.} First, we estimate the number $N_n$ of
parallelograms of the lattice $\La$ in the half-annulus $P^+(0,R_n,2
R_n)$ for $n\geq 2$. We have
\begin{equation}\label{hd5p11}
N_n \asymp \frac{4\pi R_n^2-\pi R_n^2}{2a^2(\La)}= \frac{3\pi
R_n^2}{2a^2(\La)},
\end{equation}
where $a^2(\La)$ is a measure of the parallelogram of $\La$. Recall
that in Definition \ref{konstr4} we considered the segments
$$
U(b_{l,m},\vep)=\{z\in\oc\colon -\frac{3\pi}{8}\leq \Arg
(z-b_{l,m})\leq\frac{3\pi}{8},\ |z-b_{l,m}|\leq\vep\},
$$
where $b_{l,m}\in\La$ and $\vep>0$ as in \eqref{epsilon2}. Hence,
$\vol(U(b_{l,m},\vep))=3\pi\vep^2/8$.

Fix $n\geq 2$ and $A_n\in\mathcal{A}_n$. There exist $l,m \in \mathbb{Z}$ such that $A_n =
h_n^{-1}(U(b_{l,m}^{(n)},\vep))$, where
$U(b_{l,m}^{(n)},\vep)\sbt P^+(0,R_n,2R_n)$. Moreover, for each
$A_k\in \mathcal{A}_{n+1}$  there are $l'=l'(k),m'=m'(k) \in\mathbb{Z}$ such that
$A_k=h_{n+1}^{-1}(U(b_{l',m'}^{(n+1)},\vep))$, where
$U(b_{l',m'}^{(n+1)},\vep)\sbt P^+(0,R_{n+1},2R_{n+1})$. To
simplify the formulas we denote $b_{l,m}^{(n)}$ by $b_n$. There
are finitely many sets $A_k\in \mathcal{A}_{n+1}$ contained in
$A_n$. We denote by $b_k$ the pole corresponding to $A_k$. Let
$\b_n:=h_n^{-1}(b_n)\in A_n, \b_k:=h_{n+1}^{-1}(b_k)\in A_k$. Lemma
\ref{konstr6} implies that $h_n$ are conformal on $A_n$. Using
\eqref{dystodw}, we get
$$ L(h_n,A_n)=L(h_n^{-1},h_n(A_n)).$$
Hence,
\begin{equation}\label{hd5p1}
\begin{aligned}
 \vol  (A_n) &  = \,\, \vol (h_n^{-1}(U(b_n,\vep)))
=\iint_{h_n^{-1}(U(b_n,\vep))} d\b\\
&  =  \,\,\iint_{U(b_n,\vep)} \left|(h_n^{-1})'(z)\right|^2 dz \leq
\iint_{U(b_n,\vep)} \left(\sup_{z\in U(b_n,\vep)}
|(h_n^{-1})'(z)|\right)^2 dz \\
&  =  \,\,\vol(U(b_n,\vep)) \left(L(h_n^{-1},U(b_n,\vep)) \inf_{z\in
U(b_n,\vep)}
|(h_n^{-1})'(z)|\right)^2  \\
& \leq   \frac{3\pi\vep^2}{8} \left(L(h_n,A_n)
|(h_n^{-1})'(b_n)|\right)^2  =  \frac{3\pi\vep^2}{8}
\left(\frac{L(h_n,A_n)}{|h'_n(\b_n)|}\right)^2.
\end{aligned}
\end{equation}
Set $P_{n+1}:=P^+(0,R_{n+1},2R_{n+1})$.
\begin{equation}\label{hd5p2}
\begin{aligned}
&  \vol ( \mathcal{U}_{n+1}  \cap A_n)  =   \sum_{A_k\sbt A_n} \vol (A_k) = \sum_{b_k\in P_{n+1}} \vol (h_{n+1}^{-1}(U(b_k,\vep))) \\
 & =  \sum_{b_k\in P_{n+1}} \iint_{U(b_k,\vep)} \left|(h_{n+1}^{-1})'(z)\right|^2 dz
\geq \sum_{b_k\in P_{n+1}} \iint_{U(b_k,\vep)} \left(\inf_{z\in U(b_k,\vep)} |(h_{n+1}^{-1})'(z)|\right)^2 dz  \\
  & =   \frac{3\pi\vep^2}{8} \sum_{b_k\in P_{n+1}} \left(\frac{\sup_{z\in
U(b_k,\vep)}
|(h_{n+1}^{-1})'(z)|}{L(h_{n+1}^{-1},U(b_k,\vep))}\right)^2 \geq
\frac{3\pi\vep^2}{8} \sum_{b_k\in P_{n+1}} \left(\frac{ |(h_{n+1}^{-1})'(b_k)|}{L(h_{n+1}^{-1},U(b_k,\vep))}\right)^2  \\
 & =   \frac{3\pi\vep^2}{8} \sum_{\b_k\in A_k \sbt A_n}
\left(L(h_{n+1},A_k) |h'_{n+1}(\b_k)|\right)^{-2}.
\end{aligned}
\end{equation}
Now, using \eqref{hd5p1} and \eqref{hd5p2}, we estimate the density
of the sets $\mathcal{U}_{n+1}\cap A_n$ in $A_n$.
\begin{equation}\label{hd5p3}
\begin{split}
\frac{\vol (\mathcal{U}_{n+1} \cap A_n)}{\vol (A_n)} & \geq \frac{
\sum_{\b_k\in A_k \sbt A_n} \left(L(h_{n+1},A_k)
|h'_{n+1}(\b_k)|\right)^{-2} }{\left(\frac{L(h_n,A_n)}{|h'_n(\b_n)|}\right)^2 }  \\ & =
\frac{|h'_n(\b_n)|^2}{(L(h_n,A_n))^2} \sum_{\b_k\in A_k \sbt A_n}
\left(L(h_{n+1},A_k) |h'_{n+1}(\b_k)|\right)^{-2}.
\end{split}
\end{equation}
Lemma \ref{hd1} and inequalities \eqref{hd2p3} give
\begin{equation}\label{hd5p4}
|h'_n(\b_n)| \geq \frac{R_1}{2(1+r)}\left|\prod_{j=1}^{n-1}
g_{\b_n}'(g_{\b_n}^j(c_1))\right|
\end{equation}
and
\begin{equation}\label{hd5p5}
|h'_{n+1}(\b_k)| \leq \frac{5R_1}{2(1-r)}\left|\prod_{j=1}^{n}
g_{\b_k}'(g_{\b_k}^j(c_1))\right|.
\end{equation}
It follows from Lemma \ref{hd3} that
\begin{equation}\label{hd5p6}
(L(h_n,A_n))^2 \leq \left(\frac{1+r}{1-r}\right)^2 5^2 2^{3(n-1)}
\end{equation}
and
\begin{equation}\label{hd5p7}
(L(h_{n+1},A_k))^2 \leq \left(\frac{1+r}{1-r}\right)^2 5^2 2^{3n}.
\end{equation}
Plugging \eqref{hd5p4}-\eqref{hd5p7} into \eqref{hd5p3}, we have
\begin{eqnarray}\label{hd5p8}
\begin{aligned}
& \frac{\vol (\mathcal{U}_{n+1} \cap A_n)}{\vol (A_n)} \geq \\
& \geq
\frac{\left(\frac{R_1}{2(1+r)}\right)^2\left|\prod_{j=1}^{n-1}
g_{\b_n}'(g_{\b_n}^j(c_1))\right|^2}{\left(\frac{1+r}{1-r}\right)^2
5^2 2^{3(n-1)}} \sum_{\b_k\in A_k \sbt A_n}
\frac{1}{\left(\frac{1+r}{1-r}\right)^2 5^2 2^{3n}
\left(\frac{5R_1}{2(1-r)}\right)^2\left|\prod_{j=1}^{n}
g_{\b_k}'(g_{\b_k}^j(c_1))\right|^2}  \\ &  =
\left(\frac{1-r}{1+r}\right)^6 \frac{1}{5^6 2^{3(2n-1)}}
\left|\prod_{j=1}^{n-1} g_{\b_n}'(g_{\b_n}^j(c_1))\right|^2
\sum_{\b_k\in A_k \sbt A_n} \frac{1}{\left|\prod_{j=1}^{n}
g_{\b_k}'(g_{\b_k}^j(c_1))\right|^2}  \\  & =
\left(\frac{1-r}{1+r}\right)^6 \frac{1}{5^6 2^{3(2n-1)}}
\sum_{\b_k\in A_k \sbt A_n} \frac{\prod_{j=1}^{n-1} \left|
g_{\b_n}'(g_{\b_n}^j(c_1))\right|^2}{\prod_{j=1}^{n-1} \left|
g_{\b_k}'(g_{\b_k}^j(c_1))\right|^2 } \cdot \frac{1}{\left|
g_{\b_k}'(g_{\b_k}^n(c_1))\right|^2}  \\ & =
\left(\frac{1-r}{1+r}\right)^6 \frac{1}{5^6 2^{3(2n-1)}}
\sum_{\b_k\in A_k \sbt A_n} \left(\prod_{j=1}^{n-1} \frac{\left|
g_{\b_n}'(g_{\b_n}^j(c_1))\right|}{\left|
g_{\b_k}'(g_{\b_k}^j(c_1))\right| }\right)^2 \cdot \frac{1}{\left|
g_{\b_k}'(g_{\b_k}^n(c_1))\right|^2}.
\end{aligned}
\end{eqnarray}
For each $j=1,2,\ldots, n-1$
$$ g_{\b_n}(g_{\b_n}^j(c_1))=g_{\b_n}^{j+1}(c_1)=h_{j+1}(\b_n)\in P^+(0,R_{j+1},2 R_{j+1})$$ and $$ g_{\b_k}(g_{\b_k}^j(c_1))=g_{\b_k}^{j+1}(c_1)=h_{j+1}(\b_k)\in P^+(0,R_{j+1},2 R_{j+1}),$$ since $\b_n\in A_n\sbt A_{j+1}$ and $\b_k\in A_k\sbt A_n\sbt A_{j+1}$. Thus, by \eqref{nierpochH}, for $j=1,2,\ldots n-1$ we have
$$
|g'_{\b_n}(g_{\b_n}^j(c_1))|\geq \frac{C_2 R^{3/2}_{j+1}}{C_1^{3/2}}
\ \ \textrm{ and }\ \ |g'_{\b_k}(g_{\b_k}^j(c_1))|\leq \frac{2^{3/2}
C_2 R^{3/2}_{j+1}}{C_1^{3/2}}.
$$
This implies that
\begin{equation}\label{hd5p9}
\frac{\left| g_{\b_n}'(g_{\b_n}^j(c_1))\right|}{\left|
g_{\b_k}'(g_{\b_k}^j(c_1))\right| } \geq \frac{1}{2^{3/2}},\
j=1,2,\ldots n-1.
\end{equation}
Analogously,
$$ g_{\b_k}(g_{\b_k}^n(c_1))=g_{\b_k}^{n+1}(c_1)=h_{n+1}(\b_k)\in
P^+(0,R_{n+1},2 R_{n+1})$$ as $\b_k\in A_k\in \mathcal{A}_{n+1}$. By
applying this to \eqref{nierpochH}, we get
\begin{equation}\label{hd5p10}
|g'_{\b_k}(g_{\b_k}^n(c_1))|\leq \frac{2^{3/2} C_2
R^{3/2}_{n+1}}{C_1^{3/2}}.
\end{equation}
Putting \eqref{hd5p9}, \eqref{hd5p10} into \eqref{hd5p8} and by
\eqref{hd5p11}, we obtain
\begin{equation*}
\begin{aligned}
\frac{\vol (\mathcal{U}_{n+1} \cap A_n)}{\vol A_n}&  \geq
\left(\frac{1-r}{1+r}\right)^6 \frac{1}{5^6 2^{3(2n-1)}}
\left(\frac{1}{2^{3/2}}\right)^{2(n-1)} \frac{C_1^3}{2^3 C_2^2
R_{n+1}^3 } \sum_{\b_k\in A_k \sbt A_n} 1 \\
 & =
\left(\frac{1-r}{1+r}\right)^6 \frac{2^3}{5^6 2^{9n}}
\frac{C_1^3}{C_2^2 R_{n+1}^3 } N_{n+1} \asymp
\left(\frac{1-r}{1+r}\right)^6 \frac{2^3}{5^6 2^{9n}}
\frac{C_1^3}{C_2^2 R_{n+1}^3 } R_{n+1}^2\\
&  = \frac{M}{2^{9n} R_{n+1}},
\end{aligned}
\end{equation*}
where $M = \frac{2^3 (1-r)^6 C_1^3}{5^6 (1+r)^6 C_2^2}$.

Similarly, we consider the case $n=1$. By Definition \ref{konstr4},
the set $\mathcal{A}_1$ has only one element, i.e. $A_1$ and its
Lebesgue measure $\vol (A_1)\leq \pi r^2$. The set $A_1$ contains
finitely many subsets $A_k\in \mathcal{A}_2$. As for $n\geq 2$, we
denote by $b_k$ the pole corresponding to $A_k$. Arguing as in
\eqref{hd5p2}, we get
$$
\vol (\mathcal{U}_{2} \cap A_1) \geq
 \frac{3\pi\vep^2}{8} \sum_{\b_k\in A_k \sbt A_1} \left(L(h_2,A_k) |h'_2(\b_k)|\right)^{-2}.
$$
Setting $n=1$ in bounds \eqref{hd5p5}, \eqref{hd5p7} we have
$$
|h'_2(\b_k)| \leq \frac{5R_1}{2(1-r)} |g_{\b_k}'(g_{\b_k}(c_1))|\ \
\textrm{ and }\ \ (L(h_2,A_k))^2 \leq \left(\frac{1+r}{1-r}\right)^2
5^2 2^3,
$$
which implies that
$$
\vol (\mathcal{U}_{2} \cap A_1) \geq
 \frac{3\pi\vep^2}{8} \sum_{\b_k\in A_k \sbt A_1} \frac{1}{\left(\frac{5R_1}{2(1-r)}\right)^2 |g_{\b_k}'(g_{\b_k}(c_1))|^2 \left(\frac{1+r}{1-r}\right)^2 5^2 2^3}.
$$
Analogously as in \eqref{hd5p10}, we obtain
$$
|g'_{\b_k}(g_{\b_k}(c_1))|\leq \frac{2^{3/2} C_2
R^{3/2}_2}{C_1^{3/2}}
$$
and we conclude that
\begin{equation*}
\begin{aligned}
\frac{\vol (\mathcal{U}_{2} \cap A_1)}{\vol (A_1)} & \geq
 \frac{3\pi\vep^2}{8\pi r^2} \sum_{\b_k\in A_k \sbt A_1} \frac{1}{\left(\frac{5R_1}{2(1-r)}\right)^2 \left(\frac{2^{3/2} C_2 R^{3/2}_2}{C_1^{3/2}}\right)^2 \left(\frac{1+r}{1-r}\right)^2 5^2 2^3} \\ & = \frac{3\vep^2 (1-r)^4 C_1^3}{2^7 5^4 r^2(1+r)^2 C_2^2 R_1^2} \frac{\sum_{\b_k\in A_k \sbt A_1} 1}{R_2^3} = M' \frac{N_2}{R_2^3} \asymp M' \frac{R_2^2}{R_2^3} =  \frac{M'}{R_2},
 \end{aligned}
\end{equation*}
where $M' = \frac{3\vep^2 (1-r)^4 C_1^3}{2^7 5^4 r^2(1+r)^2 C_2^2
R_1^2}$.
\qed

\

\ni By Lemma \ref{hd5}, the numbers $\Delta_n$ from
Proposition~\ref{mcmullen} are equal to
$$
\Delta_1 = \frac{M'}{R_2},\ \ \Delta_n = \frac{M}{2^{9n} R_{n+1}},\
n\geq 2.
$$

\fr Assembling the preceding lemmas, we may
now prove Theorem~\ref{Cantor}.

\

\fr{\sl Proof of Theorem~\ref{Cantor}.}  Lemma~\ref{hd5} implies
that for $a>a_0$ we obtain
\begin{eqnarray}\label{hd5p12}
\begin{aligned}
\sum_{j=1}^n & |\log\Delta_j|  =  |\log \Delta_1| + \sum_{j=2}^n
|\log \Delta_j| = \left|\log \frac{M'}{R_2}\right| + \sum_{j=2}^n
\left|\log \frac{M}{2^{9j} R_{j+1}}\right|  \\
= & \log (a R_1) - \log M' + \sum_{j=2}^n \log(2^{9j} a^j R_1) -
(n-1)\log M   \\
= &\log M - \log M' +n\log
R_1 -n\log M  + 9 \log 2 \sum_{j=2}^n j + \log a \sum_{j=1}^n j  \\
=& \log\frac{M}{M'} + n\log\frac{R_1}{M} + \frac{9(n+2)(n-1)}{2}
\log 2+ \frac{n(n+1)}{2}\log a.
\end{aligned}
\end{eqnarray}In view of  Lemma~\ref{hd4}, for $a>a_0$ we have
\begin{equation}\label{hd5p13}
\begin{aligned}  & |\log d_n|=  \left|\log \frac{4 \vep
(1+r)}{\left(\frac{C_2}{C_1^{3/2}}\right)^{n-1} a^\frac{3n(n-1)}{4}
R_1^\frac{3n-1}{2}} \right|  \\  & = (n-1)\log \frac{C_2}{C_1^{3/2}}
+ \frac{3n(n-1)}{4}\log a + \frac{3n-1}{2}\log R_1 - \log 4\vep
(1+r).
\end{aligned}
\end{equation}
The  final estimate follows from  \eqref{hd5p12} and \eqref{hd5p13}.
For $a>a_0$ we have
$$
\begin{aligned}
\dim_H (A(a)) &\geq 2 - \limsup_{n\to\infty} \frac{\log\frac{M}{M'}
+ n\log\frac{R_1}{M} +  \frac{9(n+2)(n-1)}{2} \log 2+
\frac{n(n+1)}{2}\log a}{(n-1)\log \frac{C_2}{C_1^{3/2}} +
\frac{3n(n-1)}{4}\log a + \frac{3n-1}{2}\log R_1 - \log 4\vep (1+r)}\\
 &  = 2 - \frac{\frac{1}{2}\log a + \frac{9}{2}\log
2}{\frac{3}{4}\log a} = 2 - \frac{2}{3} - \frac{6\log 2}{\log a} =
\frac{4}{3} - \frac{6\log 2}{\log a}. \end{aligned}$$
\qed

\

\fr Thus, the theorem stated in section 1  follows from
Theorem~\ref{Cantor}.

\

\begin{que} Is the Hausdorff dimension of the escaping set $\mathcal E$ equal to $4/3$ ?
\end{que}

\end{document}